\documentclass[11pt,a4paper]{article}
\usepackage[T1]{fontenc}
\usepackage[utf8]{inputenc}
\usepackage{lmodern}
\usepackage{amsmath,amssymb,amsthm,mathtools}
\usepackage{geometry,microtype,enumitem,booktabs,array,xcolor,hyperref}

\hypersetup{colorlinks=true,linkcolor=blue!45!black,citecolor=blue!45!black,urlcolor=blue!45!black}
\setlist{nosep,leftmargin=2em}
\newtheorem{theorem}{Theorem}[section]
\newtheorem{proposition}[theorem]{Proposition}
\newtheorem{corollary}[theorem]{Corollary}
\newtheorem{lemma}[theorem]{Lemma}
\theoremstyle{definition}
\newtheorem{definition}[theorem]{Definition}
\newtheorem{example}[theorem]{Example}
\theoremstyle{remark}
\newtheorem{remark}[theorem]{Remark}

\begin{document}

\title{\textbf{Split Kernel Multiplication Modules}}
\author{Violeta Leoreanu-Fotea\thanks{%
Faculty of Mathematics, Al. I. Cuza University, Bd. Carol I, No. 11, 700506
Ia\c{s}i, Romania. E-mail: \texttt{foteavioleta@gmail.com}} \and \"{U}nsal
Tekir\thanks{%
Department of Mathematics, Marmara University, Istanbul, T\"urkiye. E-mail: 
\texttt{utekir@marmara.edu.tr}} \and Dilara Erdemir\thanks{%
Department of Mathematics, Y\i ld\i z Technical University, Istanbul,
T\"urkiye. E-mail: \texttt{dilaraer@yildiz.edu.tr}}}
\date{}
\maketitle

\begin{abstract}
We introduce split kernel multiplication (SKM) modules, for which every submodule has the form $N=IM\oplus\ker f$, where $I$ is an ideal of $R$ and $f$ is an endomorphism of $M$. We study direct sums, localization, and connections with multiplication and retractability properties. For the idempotent version (ISKM), we prove that $R\oplus k^b$ is ISKM over a local ring $(R,\mathfrak m)$, with residue field $k$ and $b\geq1$, if and only if $\mathfrak m^2=0$. We also classify ISKM modules of the form $R^a\oplus k^b$ over local rings. For commutative Noetherian rings, we show that $R^2$ is SKM if and only if $R$ is an Artinian principal ideal ring. Finally, we prove that every finitely generated $R$-module is ISKM if and only if $R$ is a finite product of fields.
\end{abstract}



\section{Introduction}

An $R$-module $M$ is called a \emph{multiplication module} if every submodule of $M$ is of the form $IM$ for some ideal $I$ of $R$. This notion was introduced by Barnard~\cite{Barnard} and further studied by El-Bast and Smith~\cite{ElBastSmith}. Related classes include $S$-multiplication modules~\cite{AndersonEtAl}. Such modules describe submodules through ideals of the ring, but this approach does not apply to all modules.

A different approach uses kernels of endomorphisms. An $R$-module $M$ is \emph{co-epi-retractable} if $M/N$ embeds into $M$ for every submodule $N$~\cite{Ghorbani}. Equivalently, every submodule is the kernel of an endomorphism of $M$. The weaker notion of a \emph{co-epi-finite-retractable} (CEFR) module requires this property only for finitely generated submodules~\cite{Neishabouri}. In both cases, the endomorphisms need not be idempotent.

These two approaches lead to different classes of modules. For example, $\mathbb Z$ is a multiplication module, but $2\mathbb Z$ is not the kernel of any endomorphism of $\mathbb Z$. On the other hand, a vector space of dimension at least two is not a multiplication module, although every subspace is the kernel of an idempotent endomorphism.

This suggests combining the two representations. We introduce \emph{split kernel multiplication} (SKM) modules, in which every submodule has a decomposition $N=IM\oplus\ker f$, where $I$ is an ideal of $R$ and $f$ is an endomorphism of $M$. The sum is direct inside $N$, but $\ker f$ need not be a direct summand of $M$. We also consider \emph{idempotently split kernel multiplication} (ISKM) modules, where $f$ can be chosen idempotent. Thus every multiplication module is ISKM, every ISKM module is SKM, and every co-epi-retractable module is SKM. The definitions and their differences are discussed in Section~2.

Our main results concern the structure of these modules over local and Noetherian rings. Let $(R,\mathfrak m)$ be a local ring with residue field $k=R/\mathfrak m$. We classify ISKM modules of the form $R^a\oplus k^b$ without assuming that $R$ is Artinian. In particular, for $b\geq1$, we prove that $R\oplus k^b$ is ISKM if and only if $\mathfrak m^2=0$. We also characterize the SKM property of $R^a\oplus k^b$ when $\mathfrak m^2=0$.

We next study the SKM property of $R^2$ over Noetherian rings. Applying the SKM condition to submodules of the form $Je_1$ yields an embedding $R/J\oplus R\hookrightarrow R^2$. This gives a bound on the dimension of the socle of $R/J$ and leads to the conclusion that $R$ is an Artinian principal ideal ring. The converse follows from known results on CEFR modules over commutative K\"othe rings~\cite{Neishabouri} and the K\"othe--Cohen--Kaplansky characterization~\cite{Baghdari}. Related embedding criteria for Artinian principal ideal rings appear in~\cite{Alvarado2019,Cejudo2025}. Our result characterizes these rings through the SKM property of $R^2$.

We also prove that, over a nonzero commutative ring, $R^2$ is ISKM if and only if $R$ is a finite product of fields. Consequently, every finitely generated $R$-module is ISKM precisely for these rings. Further results concern direct sums, localization, and examples separating the classes under consideration.

Section~2 contains the definitions and basic examples. Section~3 deals with direct sums and localization. Section~4 studies rank obstructions and local ISKM modules. Section~5 establishes the Noetherian SKM criteria and related results for Artinian rings. Section~6 gives the global ISKM characterization and further consequences.

Throughout, all rings are nonzero, commutative, and have identity, and all modules are unital. The zero module is allowed. We write $k^b$ for the direct sum of $b$ copies of $k$, with $k^0=0$.

\section{Definitions and basic examples}

\begin{definition}
An $R$-module $M$ is called a \emph{split kernel multiplication module} (\emph{SKM module}) if every submodule $N\leq M$ can be written as $N=IM\oplus\ker f$, where $I$ is an ideal of $R$ and $f\in\operatorname{End}_R(M)$. The sum is internal in $N$, and $\ker f$ need not be a direct summand of $M$.
\end{definition}

\begin{definition}
An $R$-module $M$ is called an \emph{idempotently split kernel multiplication module} (\emph{ISKM module}) if every submodule $N\leq M$ admits such a decomposition with $f^2=f$.
\end{definition}

We recall the conditions defining the related classes of modules. In each case, the endomorphism may depend on the submodule.

\begin{center}
{\small
\renewcommand{\arraystretch}{1.22}
\begin{tabular}{@{}p{0.28\linewidth}p{0.64\linewidth}@{}}
\toprule
Class & Requirement \\
\midrule
Multiplication & $N=IM$ for every $N\leq M$ \\
Co-epi-retractable & $N=\ker f_N$ for every $N\leq M$ \\
CEFR & $N=\ker f_N$ for every finitely generated $N\leq M$ \\
SKM & $N=IM\oplus\ker f_N$ for every $N\leq M$ \\
ISKM & $N=IM\oplus\ker f_N$, with $f_N^2=f_N$, for every $N\leq M$ \\
\bottomrule
\end{tabular}
}
\end{center}

Here $I$ is an ideal of $R$ and $f_N\in\operatorname{End}_R(M)$. In the co-epi-retractable case, the kernel condition is equivalent to the embedding $M/N\hookrightarrow M$. Unlike CEFR, the SKM and ISKM conditions concern all submodules.

If $f^2=f$, then $M=\ker f\oplus\operatorname{Im}f$.
Hence every ISKM module is SKM. Moreover, if
$N=IM\oplus\ker f$ with $f^2=f$, then
$I\ker f=0$ and
$N\cap\operatorname{Im}f=I\operatorname{Im}f$.
Indeed, the first equality follows from
$IM\cap\ker f=0$, while the second follows
by applying $f$ to the decomposition of $N$.

\begin{proposition}
Every multiplication module is ISKM.
\end{proposition}

\begin{proof}
For $N=IM$, take $f=1_M$. Then $f^2=f$ and $\ker f=0$.
\end{proof}

\begin{example}[The Pr\"ufer group]
Let $p$ be a prime and $M=\mathbb{Z}(p^\infty)$, regarded as a $\mathbb{Z}$-module. Since $M$ is divisible, $(n\mathbb{Z})M=M$ for every $n\ne0$, while $(0)M=0$. Its proper nonzero subgroups are $G_j=\langle p^{-j}+\mathbb{Z}\rangle$, $j\geq1$, and $G_j=\ker(p^j1_M)$. Hence $M$ is SKM but not multiplication.

Moreover, $\operatorname{End}_{\mathbb{Z}}(M)\cong\mathbb{Z}_p$, which has no idempotents other than $0$ and $1$. Thus $M$ is not ISKM, since its proper nonzero subgroups are neither ideal multiples nor kernels of idempotent endomorphisms. The classical properties of the Pr\"ufer group used here are also discussed in~\cite{AndersonCorrection}.
\end{example}

\begin{proposition}
Every co-epi-retractable module is SKM.
\end{proposition}

\begin{proof}
Let $N\leq M$. Since $M$ is co-epi-retractable, there is an injection $g:M/N\hookrightarrow M$. If $\pi:M\to M/N$ is the canonical projection, then $\ker(g\pi)=N$. Taking $I=0$ gives the required decomposition.
\end{proof}

\begin{example} 
Let $k$ be a field and $R=k[X,Y]/(X^2,XY,Y^2)$. Write $x,y$ for the residue classes of $X,Y$, and put $\mathfrak m=(x,y)$, $M=R\oplus k$, and $N=(x)\oplus k$.

For $I=(x)$ and $\pi_1(r,c)=(r,0)$, we have $IM=(x)\oplus0$ and $\ker\pi_1=0\oplus k$. Thus $N=IM\oplus\ker\pi_1$, with both summands nonzero.

However, $N\ne JM$ for every ideal $J$ of $R$. Indeed, if $J\subseteq\mathfrak m$, then $Jk=0$, whereas $N$ contains $0\oplus k$. If $J=R$, then $JM=M\ne N$.

Also, $N$ is not the kernel of an endomorphism of $M$. Otherwise, $M/N\cong R/(x)$ would embed into $M$. The element $1+(x)$ is annihilated by $x$ but not by $y$. Yet, for $(r,c)\in M$, the equality $x(r,c)=0$ implies $r\in\mathfrak m$, and hence $y(r,c)=0$, since $\mathfrak m^2=0$. This rules out such an embedding.

Therefore $N$ admits an SKM decomposition, although it is neither an ideal multiple nor an endomorphism kernel.
\end{example}

\begin{proposition}
Every semisimple module is ISKM. In particular, every vector space over a field is ISKM.
\end{proposition}

\begin{proof}
For $N\leq M$, write $M=N\oplus K$. The projection onto $K$ along $N$ is idempotent and has kernel $N$. Taking $I=0$ gives the result.
\end{proof}

\begin{proposition}
Suppose every nonzero endomorphism of $M$ is injective. Then $M$ is SKM if and only if it is a multiplication module.
\end{proposition}

\begin{proof}
Every multiplication module is SKM. Conversely, suppose $M$ is SKM and let $N\lneq M$. Write $N=IM\oplus\ker f$. Since $N\ne M$, we have $f\ne0$. Hence $f$ is injective, so $N=IM$. The equality $M=RM$ completes the proof.
\end{proof}

\begin{lemma} 
Let $(R,\mathfrak m)$ be a local ring and $M$ an SKM module. Then every cyclic submodule $N$ of $M$ is either an ideal multiple of $M$ or the kernel of an endomorphism of $M$. If $M$ is ISKM, the second alternative can be replaced by the assertion that $N$ is a direct summand of $M$.
\end{lemma}

\begin{proof}
The case $N=0$ is clear. Otherwise, $N\cong R/J$ for some proper ideal $J$. Since $R/J$ is local, its only idempotents are $0$ and $1$. As $\operatorname{End}_R(N)\cong R/J$, the module $N$ is indecomposable.

Now write $N=IM\oplus\ker f$. Since $N$ is indecomposable, either $IM=0$ or $\ker f=0$. Thus $N=\ker f$ or $N=IM$. If $M$ is ISKM, we may choose $f^2=f$, so $M=\ker f\oplus\operatorname{Im}f$.
\end{proof}


\section{Componentwise constructions and localization}

\begin{theorem}\label{thm:comax}
Let $M=M_1\oplus M_2$, where
$\operatorname{Ann}_R(M_1)+\operatorname{Ann}_R(M_2)=R$.
Then $M$ is SKM if and only if both $M_1$ and $M_2$ are SKM.
\end{theorem}

\begin{proof}
Choose $a_i\in\operatorname{Ann}_R(M_i)$ with $a_1+a_2=1$.
Multiplication by $a_2$ and $a_1$ gives the projections onto
$M_1$ and $M_2$, respectively. Thus every submodule of $M$
has the form $N=N_1\oplus N_2$.

Moreover, $\operatorname{Hom}_R(M_1,M_2)=
\operatorname{Hom}_R(M_2,M_1)=0$.
Indeed, if $h:M_1\to M_2$, then
$h(m)=h(a_2m)=a_2h(m)=0$.
Consequently, every endomorphism of $M$ is of the form
$f=f_1\oplus f_2$.

Suppose that both $M_i$ are SKM. For
$N_i\leq M_i$, write
$N_i=I_iM_i\oplus\ker f_i$.
Set $I=a_2I_1+a_1I_2$ and $f=f_1\oplus f_2$.
Then $IM=I_1M_1\oplus I_2M_2$ and
$\ker f=\ker f_1\oplus\ker f_2$.
Hence $N=IM\oplus\ker f$.

Conversely, suppose $M$ is SKM and let $N_1\leq M_1$.
Write $N_1\oplus0=IM\oplus\ker f$.
Since $f=f_1\oplus f_2$, taking the first coordinate gives
$N_1=IM_1\oplus\ker f_1$.
Thus $M_1$ is SKM, and similarly $M_2$ is SKM.
\end{proof}

\begin{corollary}
Let $M=\bigoplus_{i=1}^tM_i$ and suppose that
$\operatorname{Ann}_R(M_i)$ are pairwise comaximal.
Then $M$ is SKM if and only if every $M_i$ is SKM.
\end{corollary}

\begin{proof}
Apply Theorem~\ref{thm:comax} inductively, using
$\operatorname{Ann}_R(\bigoplus_{i<t}M_i)
=\bigcap_{i<t}\operatorname{Ann}_R(M_i)$
and the fact that a finite intersection of ideals
comaximal with a fixed ideal remains comaximal with it.
\end{proof}

\begin{remark}\label{rem:iskmprod}
Theorem~\ref{thm:comax} and its corollary
remain valid with SKM replaced by ISKM.
Indeed, the endomorphisms in the proof
decompose componentwise, and
$f=\bigoplus_i f_i$ is idempotent if and
only if each $f_i$ is idempotent.
\end{remark}

\begin{example}
Let $p$ be a prime, $R=\mathbb Z\times\mathbb Z$, and
$M=\mathbb Z\oplus\mathbb Z(p^\infty)$, with componentwise
action. The annihilators of the two components are
$0\times\mathbb Z$ and $\mathbb Z\times0$, which are
comaximal. Since both components are SKM,
Theorem~\ref{thm:comax} shows that $M$ is SKM.

Consider $N=2\mathbb Z\oplus G_1$, where $G_1$ is the
subgroup of order $p$ in $\mathbb Z(p^\infty)$.
Taking $I=2\mathbb Z\times0$ and
$f(a,b)=(a,pb)$, we obtain
$N=IM\oplus\ker f$.

However, $N$ is not an ideal multiple of $M$, since
the second component $G_1$ is not an ideal multiple
of $\mathbb Z(p^\infty)$. It is also not the kernel
of an endomorphism of $M$, since no endomorphism
of $\mathbb Z$ has kernel $2\mathbb Z$.
\end{example}

\begin{remark}
Comaximality ensures that every submodule splits
componentwise and that there are no nonzero
homomorphisms between distinct components.
The result need not hold for arbitrary direct sums.
\end{remark}

\subsection{Localization}

\begin{theorem}\label{thm:loc}
Let $S$ be a multiplicatively closed subset of $R$ with $0\notin S$. If $M$ is SKM (respectively, ISKM), then $S^{-1}M$ is SKM (respectively, ISKM) over $S^{-1}R$.
\end{theorem}

\begin{proof}
Let $U\leq S^{-1}M$ and put $N=\{m\in M:m/1\in U\}$. Then $U=S^{-1}N$.

Suppose $N=IM\oplus\ker f$. Since localization is exact, we obtain
\[
U=(S^{-1}I)(S^{-1}M)\oplus\ker(S^{-1}f).
\]
Thus $S^{-1}M$ is SKM. If $M$ is ISKM, we may choose $f^2=f$, and then $(S^{-1}f)^2=S^{-1}f$.
\end{proof}

\begin{example}[The converse to localization fails]
Let $R=\mathbb C[X,Y]$, $M=(x,y)$, where $x=X$ and $y=Y$, and $S=\{1,x,x^2,\ldots\}$.

We first show that every nonzero endomorphism of $M$ is injective. Since $M$ is a torsion-free module of rank one, each endomorphism is multiplication by some $q\in\mathbb C(X,Y)$. Write $q=a/b$, where $a,b\in R$ are coprime. As $qx,qy\in M\subseteq R$, we have $b\mid ax$ and $b\mid ay$. Hence $b$ divides both $x$ and $y$, so $b$ is a unit. Therefore $q\in R$, and every nonzero endomorphism of $M$ is injective.

Next, $M$ is not a multiplication module. Indeed, if $(x)=IM$ for some ideal $I$, then $Iy\subseteq(x)$. Since $x$ and $y$ are coprime, $I\subseteq(x)$, giving $IM\subseteq xM=(x^2,xy)\subsetneq(x)$, a contradiction.

By the proposition in Section~2, $M$ is not SKM. However, $S^{-1}M=S^{-1}R$, since $x$ is invertible in $S^{-1}R$. Thus $S^{-1}M$ is a multiplication module and hence SKM.
\end{example}

\begin{corollary}
Let $R=R_1\times\cdots\times R_t$ be a finite product of local rings. Then an $R$-module $M$ is SKM if and only if $M_{\mathfrak m}$ is SKM over $R_{\mathfrak m}$ for every maximal ideal $\mathfrak m$ of $R$.
\end{corollary}

\begin{proof}
Let $e_1,\ldots,e_t$ be the standard central idempotents of $R$. Then $M=\bigoplus_{i=1}^t e_iM$, and each component is naturally a module over $R_i$. Since the factors are local, the localizations of $M$ at the maximal ideals of $R$ are precisely the modules $e_iM$. The result follows from the finite direct-sum corollary.
\end{proof}


\section{Obstructions and ISKM classification over local rings}

\begin{proposition}
If $R$ contains a regular nonunit, then $R^n$ is not SKM for any $n\geq2$.
\end{proposition}

\begin{proof}
Let $a$ be a regular nonunit, $M=R^n$, and $N=Ra e_1$. Suppose $N=IM\oplus\ker f$. Since $IM\subseteq N$, the second coordinate gives $I=0$. Thus $N=\ker f$, so $M/N$ embeds into $M$. However, $e_1+N$ is nonzero and annihilated by $a$, whereas multiplication by $a$ is injective on $M$, a contradiction.
\end{proof}

\begin{theorem} 
\label{thm:rank}
Let $D$ be an integral domain that is not a field. Then no torsion-free $D$-module of rank at least two is SKM.
\end{theorem}

\begin{proof}
Choose linearly independent elements $x,y\in M$ and a nonzero nonunit $a\in D$. Put $N=Dax$. If $N=IM\oplus\ker f$, then $Iy\subseteq Dax$ implies $I=0$. Hence $N=\ker f$, and $M/N$ embeds into $M$.

But $x+N$ is nonzero and annihilated by $a$. This contradicts the torsion-freeness of $M$.
\end{proof}

\begin{remark}
Consequently, for $n\geq2$, the module $D^n$ is SKM if and only if $D$ is a field. We next consider local rings whose maximal ideals have square zero.
\end{remark}

\begin{theorem} 
\label{thm:construct}
Let $(R,\mathfrak m)$ be a local ring with $\mathfrak m^2=0$ and residue field $k=R/\mathfrak m$. Then $R\oplus k^b$ is ISKM for every $b\geq0$.
\end{theorem}

\begin{proof}
Let $M=R\oplus k^b$, $N\leq M$, and let $\pi:M\to R$ be the first projection.

\emph{Case 1: $\pi(N)$ contains a unit.}
Then $(1,c_0)\in N$ for some $c_0\in k^b$. Put $V=\{c\in k^b:(0,c)\in N\}$. Every $(r,c)\in N$ satisfies $(0,c-\bar r c_0)\in N$, where $\bar r$ is the image of $r$ in $k$. Hence $N=R(1,c_0)\oplus(0,V)$.

Choose $k^b=V\oplus W$. Then $M=N\oplus(0,W)$. The projection onto $(0,W)$ along $N$ is idempotent with kernel $N$, giving $N=0M\oplus\ker f$.

\emph{Case 2: $\pi(N)\subseteq\mathfrak m$.}
Now $N\subseteq\mathfrak m\oplus k^b$, which is a $k$-vector space. Define
\[
J=\{r\in\mathfrak m:(r,0)\in N\},
\qquad
V=\{c\in k^b:(r,c)\in N
\text{ for some }r\in\mathfrak m\}.
\]
Since $\mathfrak m^2=0$, the subspace $J$ is an ideal of $R$.

Choose a $k$-linear section of $N\to V$. It determines a linear map $\theta_0:V\to\mathfrak m$ such that
\[
N=(J\oplus0)\oplus
\{(\theta_0(c),c):c\in V\}.
\]

Choose $k^b=V\oplus W$, extend $\theta_0$ to $\theta:k^b\to\mathfrak m$ by setting $\theta|_W=0$, and let $P$ be the projection onto $W$ along $V$.

Define $f:M\to M$ by $f(r,c)=(r-\theta(c),P(c))$. The map $f$ is $R$-linear, since $\theta$ is
$k$-linear and both $\mathfrak m$ and $k^b$
are annihilated by $\mathfrak m$.
Moreover, $P^2=P$ and $\theta P=0$,
so $f^2=f$.

Since $\ker P=V$, we have
\[
\ker f=\{(\theta_0(c),c):c\in V\}.
\]
Finally, $JM=J\oplus0$, and therefore $N=JM\oplus\ker f$.
\end{proof}

\begin{remark}
The module $R\oplus k$ from the mixed decomposition example is therefore ISKM when $\mathfrak m^2=0$. We now prove that this condition is also necessary when $b\geq1$.
\end{remark}

\begin{theorem} 
\label{thm:sharp}
Let $(R,\mathfrak m)$ be a local ring, $k=R/\mathfrak m$, and $b\geq1$. Then $R\oplus k^b$ is ISKM if and only if $\mathfrak m^2=0$.
\end{theorem}

\begin{proof}
Sufficiency follows from Theorem~\ref{thm:construct}.

Conversely, suppose $M=R\oplus k^b$ is ISKM. Let $u\in\mathfrak m$ and choose $0\ne e\in k^b$. Put $z=(u,e)$ and $N=Rz$. Since $u\in\mathfrak m$, the submodule $N$ is proper.

We first show that $N$ is not an ideal multiple of $M$. If $I\subseteq\mathfrak m$, then $IM=I\oplus0$, whereas $N$ has a nonzero second component. If $I=R$, then $IM=M\ne N$.

By the cyclic submodule dichotomy, $N$ is a direct summand of $M$. Hence there is a retraction $q:M\to N$ with $q(z)=z$. Write $q(1,0)=cz$ and $q(0,e)=dz$ for some $c,d\in R$. Since $z=u(1,0)+(0,e)$, we obtain $(uc+d-1)z=0$.

As the second component of $z$ is nonzero, $\operatorname{Ann}_R(z)\subseteq\mathfrak m$. Thus $uc+d-1\in\mathfrak m$. Since $u\in\mathfrak m$, it follows that $d-1\in\mathfrak m$, so $d$ is a unit.

Now let $t\in\mathfrak m$. Since $t(0,e)=0$, we have
\[
0=q(t(0,e))=tq(0,e)=tdz.
\]
As $d$ is invertible, $tz=0$, and the first component gives $tu=0$. Since $t,u\in\mathfrak m$ were arbitrary, $\mathfrak m^2=0$.
\end{proof}

\begin{proposition} 
\label{prop:nak}
Let $(R,\mathfrak m)$ be a local ring and $M$ a finitely generated ISKM module. Then every submodule $N\subseteq\mathfrak m M$ is of the form $N=IM$ for some ideal $I$ of $R$.
\end{proposition}

\begin{proof}
Write $N=IM\oplus K$, where $K=\ker f$ for some idempotent endomorphism $f$ of $M$. Since $K$ is a direct summand of the finitely generated module $M$, it is finitely generated.

Let $p:M\to K$ be the corresponding projection. As $K\subseteq\mathfrak m M$, we have
\[
K=p(K)\subseteq p(\mathfrak m M)=\mathfrak m K.
\]
Thus $K=\mathfrak m K$, and Nakayama's lemma gives $K=0$. Consequently, $N=IM$.
\end{proof}

\begin{theorem} 
\label{thm:localiskm}
Let $(R,\mathfrak m)$ be a local ring, $k=R/\mathfrak m$, and $a,b\geq0$. Then $M=R^a\oplus k^b$ is ISKM if and only if one of the following holds:
\begin{enumerate}[label=\textup{(\roman*)}]
\item $R$ is a field;
\item $a=0$;
\item $a=1$ and $b=0$;
\item $a=1$, $b\geq1$, and $\mathfrak m^2=0$.
\end{enumerate}
In particular, if $R$ is not a field and $a\geq2$, then $M$ is not ISKM.
\end{theorem}

\begin{proof}
If $R$ is a field or $a=0$, then $M$ is semisimple and hence ISKM. If $a=1$ and $b=0$, then $M=R$ is multiplication. Case (iv) follows from Theorem~\ref{thm:sharp}.

Conversely, suppose $R$ is not a field and $a\geq2$. Choose $0\ne c\in\mathfrak m$ and let $N=Rce_1\subseteq\mathfrak m M$. If $M$ were ISKM, Proposition~\ref{prop:nak} would give $N=IM$ for some ideal $I$.

For every $r\in I$, we have $re_2\in IM=N$. Since all elements of $N$ have zero second free coordinate, $r=0$. Hence $I=0$, contradicting $N\ne0$.

The remaining case $a=1$, $b\geq1$ is characterized by Theorem~\ref{thm:sharp}.
\end{proof}


\section{Socle methods and SKM over Noetherian rings}

For an $R$-module $M$, the socle $\operatorname{Soc}(M)$ is the sum of its simple submodules. If $(R,\mathfrak m)$ is local with residue field $k=R/\mathfrak m$, then $\operatorname{Soc}(M)=(0:_M\mathfrak m)$, which is naturally a $k$-vector space. In particular, $\operatorname{Soc}(R)=\operatorname{Ann}_R(\mathfrak m)$.

\subsection{Socle obstructions and quotient bounds}

\begin{theorem} 
\label{thm:socle}
Let $(R,\mathfrak m)$ be a local ring with residue field $k$. If $\dim_k\operatorname{Soc}(R)\geq2$, then $R^a\oplus k^b$ is not SKM for any $a\geq2$ and $b\geq0$.
\end{theorem}

\begin{proof}
Choose linearly independent elements $x,y\in\operatorname{Soc}(R)$. Since $x\ne0$ and $\mathfrak m x=0$, we have $\operatorname{Ann}_R(x)=\mathfrak m$. Moreover, $y\notin Rx$.

Suppose $M=R^a\oplus k^b$ is SKM and let $N=Rxe_1$, where $e_1,e_2$ are the first two free basis vectors. Write $N=IM\oplus\ker f$.

Since $IM\subseteq N$, we have $re_2\in N$ for every $r\in I$, which implies $I=0$. Hence $N=\ker f$, and $M/N$ embeds into $M.$ In particular, $R/(x)$ embeds into $M$.

Let $h:R/(x)\to M$ be such an embedding and write $h(1+(x))=(r_1,\ldots,r_a,c)$, with $c\in k^b$. Since $xh(1+(x))=0$, we obtain $r_i\in\mathfrak m$ for every $i$.

But $y\mathfrak m=0$ and $yk^b=0$, so $yh(1+(x))=0$. This contradicts injectivity, since $y+(x)\ne0$. Therefore $M$ is not SKM.
\end{proof}

\begin{example}
\label{ex:nonsquare}
Let $k$ be a field and $R=k[X,Y]/(X^2,XY,Y^3)$, with $\mathfrak m=(x,y)$.

The elements $1,x,y,y^2$ form a $k$-basis of $R$, and $\mathfrak m^2=(y^2)\ne0$. However, $\operatorname{Soc}(R)=(x,y^2)$ has dimension two over $k$.

By Theorem~\ref{thm:socle}, the module $R^a\oplus k^b$ is not SKM whenever $a\geq2$ and $b\geq0$. Thus the socle obstruction also applies when $\mathfrak m^2\ne0$.
\end{example}

\begin{lemma}
\label{lem:socbudget}
Let $(R,\mathfrak m)$ be a local ring and suppose $s=\dim_k\operatorname{Soc}(R)<\infty$, where $k=R/\mathfrak m$. If $M=R^a\oplus k^b$ is SKM with $a\geq2$ and $b\geq0$, then every ideal $J$ of $R$ satisfies
\[
\dim_k\operatorname{Soc}(R/J)\leq s.
\]
\end{lemma}

\begin{proof}
Let $e_1,e_2$ be distinct free basis vectors and set $N=Je_1$. In an SKM decomposition $N=IM\oplus\ker f$, the inclusion $IM\subseteq N$ implies $I=0$, since $re_2\in N$ for every $r\in I$.

Thus $N=\ker f$, and $M/N$ embeds into $M$. Since $M/N\cong R/J\oplus R^{a-1}\oplus k^b$, we obtain an embedding
\[
R/J\oplus R^{a-1}\oplus k^b
\hookrightarrow R^a\oplus k^b.
\]

An injective homomorphism maps elements annihilated by $\mathfrak m$ injectively into elements annihilated by $\mathfrak m$. Hence it induces an injection
\[
\operatorname{Soc}(R/J)\oplus
\operatorname{Soc}(R)^{a-1}\oplus k^b
\hookrightarrow
\operatorname{Soc}(R)^a\oplus k^b.
\]

The target has dimension $as+b$ over $k$. Comparing dimensions gives
$\dim_k\operatorname{Soc}(R/J)+(a-1)s+b\leq as+b$,
and therefore $\dim_k\operatorname{Soc}(R/J)\leq s$.
\end{proof}

\begin{example}
\label{ex:onesocle}
Let $R=k[X,Y]/(X^2,Y^2)$, with maximal ideal
$\mathfrak m=(x,y)$. The elements $1,x,y,xy$ form
a $k$-basis of $R$, and
$\operatorname{Soc}(R)=(xy)$ has dimension one,
while $\mathfrak m^2=(xy)\ne0$.

However, $R/\mathfrak m^2$ has a two-dimensional
socle, generated by the images of $x$ and $y$.
By Lemma~\ref{lem:socbudget}, $R^a\oplus k^b$
is not SKM for any $a\geq2$ and $b\geq0$.

Thus a one-dimensional socle does not guarantee
the SKM property.
\end{example}

The preceding results provide necessary conditions
for the SKM property. For sufficiency, we use
two known results. Neishabouri, Tolooei and
Bagheri~\cite{Neishabouri} proved that every
module over a commutative K\"othe ring is CEFR.
Moreover, by~\cite[Theorem~1.2]{Baghdari},
a commutative ring is K\"othe if and only if
it is an Artinian principal ideal ring.

\begin{corollary}
\label{cor:pir}
Let $R$ be a commutative Artinian principal ideal ring.
Then every finitely generated $R$-module is SKM.
In fact, every submodule is the kernel of an
endomorphism.
\end{corollary}

\begin{proof}
Every finitely generated module over an Artinian
ring is Noetherian, so all its submodules are
finitely generated. Since $R$ is a K\"othe ring,
the CEFR result of~\cite{Neishabouri} applies.
Thus every submodule is an endomorphism kernel,
and the SKM condition follows by taking $I=0$.
\end{proof}


\subsection{A Noetherian rank-two characterization}

\begin{theorem}
\label{thm:noethlocal}
Let $(R,\mathfrak m)$ be a nonzero Noetherian local ring, with residue field $k=R/\mathfrak m$, and let $a\geq2$, $b\geq0$. Then $R^a\oplus k^b$ is SKM if and only if $R$ is an Artinian principal ideal ring.
\end{theorem}

\begin{proof}
Suppose $M=R^a\oplus k^b$ is SKM. Since $R$ is Noetherian, its socle has finite dimension $s$ over $k$.

Applying Lemma~\ref{lem:socbudget} to $J=\mathfrak m$ gives $s\geq1$. By Theorem~\ref{thm:socle}, we also have $s\leq1$. Hence $s=1$.

Next, take $J=\mathfrak m^2$. Since $\mathfrak m/\mathfrak m^2$ is contained in $\operatorname{Soc}(R/\mathfrak m^2)$, Lemma~\ref{lem:socbudget} yields
\[
\dim_k(\mathfrak m/\mathfrak m^2)\leq1.
\]

If $\mathfrak m/\mathfrak m^2=0$, Nakayama's lemma gives $\mathfrak m=0$, so $R$ is a field.

Otherwise, choose $t\in\mathfrak m$ whose image generates $\mathfrak m/\mathfrak m^2$. Then $\mathfrak m=(t)+\mathfrak m^2$. Applying Nakayama's lemma to $\mathfrak m/(t)$ gives $\mathfrak m=(t)$.

Since $s=1$, choose $0\ne u\in\operatorname{Soc}(R)$. Then $tu=0$. By the Krull intersection theorem,
$\bigcap_{j\geq0}(t^j)=0$.
Thus there exists $j\geq0$ such that $u\in(t^j)\setminus(t^{j+1})$.

Write $u=t^jv$. If $v\in\mathfrak m=(t)$, then $u\in(t^{j+1})$, a contradiction. Hence $v$ is a unit. From $tu=0$, we obtain $t^{j+1}v=0$, and therefore $t^{j+1}=0$.

Thus $\mathfrak m$ is nilpotent. Since $R$ is Noetherian and local, it follows that $R$ is Artinian.

We now show that every ideal is principal. Let $J\ne0$ be an ideal, and choose the largest integer $i$ such that $J\subseteq(t^i)$. There exists $w\in J\setminus(t^{i+1})$. Writing $w=t^iv$, we see that $v$ is a unit. Consequently, $(t^i)\subseteq J\subseteq(t^i)$, so $J=(t^i)$. Hence $R$ is an Artinian principal ideal ring.

Conversely, if $R$ is an Artinian principal ideal ring, Corollary~\ref{cor:pir} shows that every finitely generated $R$-module is SKM. In particular, $M$ is SKM.
\end{proof}


\subsection{The global Noetherian characterization}

\begin{theorem}
\label{thm:noethglobal}
Let $R$ be a nonzero commutative Noetherian ring. The following conditions are equivalent:
\begin{enumerate}[label=\textup{(\roman*)}]
\item Every finitely generated $R$-module is SKM;
\item $R^2$ is SKM;
\item $R$ is an Artinian principal ideal ring.
\end{enumerate}
\end{theorem}

\begin{proof}
Clearly, (i) implies (ii).

Suppose $R^2$ is SKM. For every maximal ideal $\mathfrak m$ of $R$, Theorem~\ref{thm:loc} shows that $(R_{\mathfrak m})^2$ is SKM. By Theorem~\ref{thm:noethlocal}, $R_{\mathfrak m}$ is an Artinian principal ideal ring.

Let $\mathfrak p$ be a prime ideal of $R$ and choose a maximal ideal $\mathfrak m$ containing $\mathfrak p$. Since $R_{\mathfrak m}$ is Artinian local, its only prime ideal is $\mathfrak mR_{\mathfrak m}$. Hence $\mathfrak pR_{\mathfrak m}=\mathfrak mR_{\mathfrak m}$, and the prime ideal correspondence gives $\mathfrak p=\mathfrak m$.

Thus every prime ideal of $R$ is maximal. Since $R$ is Noetherian, it is Artinian and decomposes into a finite product of Artinian local rings. Each factor is isomorphic to some $R_{\mathfrak m}$ and is therefore a principal ideal ring. Consequently, $R$ is an Artinian principal ideal ring.

Finally, (iii) implies (i) by Corollary~\ref{cor:pir}.
\end{proof}

\begin{corollary}
\label{thm:artlocal}
Let $(R,\mathfrak m)$ be a nonzero Artinian local ring, $k=R/\mathfrak m$, and $a\geq2$, $b\geq0$. Then $R^a\oplus k^b$ is SKM if and only if $R$ is a principal ideal ring.
\end{corollary}

\begin{proof}
This follows from Theorem~\ref{thm:noethlocal}, since every Artinian ring is Noetherian.
\end{proof}

\begin{corollary}
\label{thm:artglobal}
Let $R$ be a nonzero commutative Artinian ring. Then the following are equivalent:
\begin{enumerate}[label=\textup{(\roman*)}]
\item Every finitely generated $R$-module is SKM;
\item $R^2$ is SKM;
\item $R$ is a principal ideal ring.
\end{enumerate}
\end{corollary}

\begin{proof}
Apply Theorem~\ref{thm:noethglobal}.
\end{proof}

\begin{remark}
\label{rem:artnovelty}
Related characterizations of Artinian principal
ideal rings appear in~\cite{Alvarado2019,Cejudo2025},
under different module-theoretic conditions.
Theorem~\ref{thm:noethglobal} provides a
characterization in terms of the SKM property
of $R^2$.
\end{remark}

\begin{remark}
The classification for $a\geq2$ differs from the local ISKM results for $a\leq1$. Indeed, $k^b$ is semisimple, $R$ is multiplication, and $R\oplus k^b$, with $b\geq1$, is ISKM precisely when $\mathfrak m^2=0$.

For example, $R=\mathbb Z/8\mathbb Z$ is an Artinian principal ideal ring. Thus $R^2$ is SKM, but it is not ISKM by Theorem~\ref{thm:localiskm}.
\end{remark}

\subsection{Square-zero consequences}

We conclude with two consequences of the preceding results for local rings with square-zero maximal ideals.

\begin{corollary}
\label{thm:completeSKM}
Let $(R,\mathfrak m)$ be a local ring with $\mathfrak m^2=0$ and residue field $k=R/\mathfrak m$. For all $a,b\geq0$, the module $R^a\oplus k^b$ is SKM if and only if $a\leq1$ or $\dim_k\mathfrak m\leq1$.
\end{corollary}

\begin{proof}
If $a=0$, the module is semisimple. If $a=1$, it is ISKM by Theorem~\ref{thm:construct}.

Suppose $a\geq2$. If $\dim_k\mathfrak m\geq2$, then $\operatorname{Soc}(R)=\mathfrak m$, since $\mathfrak m^2=0$ and $\mathfrak m\ne0$. Theorem~\ref{thm:socle} shows that $R^a\oplus k^b$ is not SKM.

Conversely, if $\dim_k\mathfrak m\leq1$, then $R$ is either a field or a local Artinian principal ideal ring with ideals $0,\mathfrak m,R$. Hence Corollary~\ref{cor:pir} gives the SKM property.
\end{proof}

\begin{corollary}
\label{cor:R2skm}
Let $(R,\mathfrak m)$ be a local ring with $\mathfrak m^2=0$ and residue field $k=R/\mathfrak m$. The following are equivalent:
\begin{enumerate}[label=\textup{(\roman*)}]
\item Every finitely generated $R$-module is SKM;
\item $R^2$ is SKM;
\item $\dim_k\mathfrak m\leq1$.
\end{enumerate}
\end{corollary}

\begin{proof}
Clearly, (i) implies (ii). By Corollary~\ref{thm:completeSKM}, (ii) implies (iii). Finally, (iii) implies that $R$ is an Artinian principal ideal ring, so Corollary~\ref{cor:pir} gives (i).
\end{proof}


\subsection{Examples separating the module classes}

We give an example distinguishing ISKM from multiplication and CEFR.

\begin{proposition}
\label{prop:noncefr}
Let $(R,\mathfrak m)$ be a local ring with
$\mathfrak m^2=0$, residue field $k=R/\mathfrak m$,
and $\dim_k\mathfrak m\geq2$.
Then $R\oplus k^b$ is ISKM for every $b\geq1$,
but is neither multiplication nor CEFR.
\end{proposition}

\begin{proof}
Put $M=R\oplus k^b$. By Theorem~\ref{thm:construct},
$M$ is ISKM.

Let $0\ne V\leq k^b$. If $0\oplus V=IM$ for some
ideal $I$, the first coordinate implies $I=0$,
a contradiction. Thus $M$ is not multiplication.

Choose linearly independent $x,y\in\mathfrak m$
and put $N=Rx\oplus0$. Suppose $N=\ker f$ for
some endomorphism $f$ of $M$. Then $M/N$ embeds
into $M$, and hence so does $R/(x)$.

Let $h:R/(x)\to M$ be an embedding and write
$h(1+(x))=(r,c)$. Since $xh(1+(x))=0$,
we have $xr=0$. As $\mathfrak m^2=0$ and
$x\ne0$, this implies $r\in\mathfrak m$.

Consequently, $y(r,c)=0$, whereas
$y+(x)\ne0$. This contradicts the injectivity
of $h$. Thus $N$ is not an endomorphism kernel.
Since $N$ is cyclic, $M$ is not CEFR.
\end{proof}

\begin{example}
Let $R=\mathbb F_2[X,Y]/(X^2,XY,Y^2)$,
with $\mathfrak m=(x,y)$ and
$k=R/\mathfrak m\cong\mathbb F_2$.

The ring $R$ has eight elements,
$\mathfrak m^2=0$, and $\dim_k\mathfrak m=2$.
Hence $M=R\oplus k$ has 16 elements and
is ISKM, but is neither multiplication
nor CEFR, by Proposition~\ref{prop:noncefr}.

Moreover, $R^2$ is not SKM by
Corollary~\ref{thm:completeSKM}.
\end{example}


\section{Global ISKM rigidity and concluding comparison}

We now distinguish SKM from ISKM and establish a global characterization of ISKM. We also record some consequences for products and free modules.

\begin{example}
Let $R=\mathbb Z/8\mathbb Z$ and $M=R^2$. Since $R$ is an Artinian principal ideal ring, Corollary~\ref{cor:pir} shows that $M$ is SKM.

For $N=2Re_1$, the endomorphism $f(u,v)=(4u,v)$ satisfies $\ker f=N$, although $f$ is not idempotent.

However, $M$ is not ISKM. Indeed, $N\subseteq2M$, and $N\ne IM$ for any ideal $I$, since the second free coordinate implies $I=0$. This contradicts Proposition~\ref{prop:nak} if $M$ is assumed ISKM.
\end{example}

\begin{theorem}
\label{thm:global}
Let $R$ be a nonzero commutative ring. The following conditions are equivalent:
\begin{enumerate}[label=\textup{(\roman*)}]
\item Every finitely generated $R$-module is ISKM;
\item $R^2$ is ISKM;
\item $R$ is a finite direct product of fields.
\end{enumerate}
\end{theorem}

\begin{proof}
Clearly, (i) implies (ii).

Suppose $R^2$ is ISKM. For any ideal $J$ of $R$, put $N=Je_1$. Write $N=IR^2\oplus\ker f$, where $f^2=f$.

Since $IR^2\subseteq Je_1$, the second coordinate implies $I=0$. Thus $N=\ker f$. As $f$ is idempotent, $R^2/N\cong\operatorname{Im}f$ is projective.

But $R^2/N\cong(R/J)\oplus R$, so $R/J$ is projective. Hence the exact sequence
\[
0\longrightarrow J\longrightarrow R
\longrightarrow R/J\longrightarrow0
\]
splits. Therefore every ideal of $R$ is a direct summand of $R$.

It follows that the regular module $R$ is semisimple. Since $R$ is cyclic, it is a finite direct sum of simple modules. Consequently, $R$ is a finite product of fields.

Conversely, if $R$ is a finite product of fields, every $R$-module is semisimple and hence ISKM.
\end{proof}


\begin{corollary}
\label{cor:free}
Let $R$ be a nonzero commutative ring and $n\geq1$.
\begin{enumerate}[label=\textup{(\alph*)}]
\item If $R$ is local, then $R^n$ is ISKM if and only if $n=1$ or $R$ is a field.
\item If $R$ is a finite product of local rings, then $R^n$ is ISKM if and only if $n=1$ or $R$ is a finite product of fields.
\end{enumerate}
\end{corollary}

\begin{proof}
Part (a) follows from Theorem~\ref{thm:localiskm}. Part (b) follows by applying Remark~\ref{rem:iskmprod} to the local factors.
\end{proof}

\begin{remark}
The following table summarizes the distinctions established above.

{\small
\begin{tabular}{@{}p{0.41\linewidth}ccc@{}}
\toprule
Module & Multiplication & SKM & ISKM\\
\midrule
$\mathbb Z(p^\infty)$ over $\mathbb Z$
& No & Yes & No\\
$R\oplus k$, $\mathfrak m^2=0$,
$\dim_k\mathfrak m\geq2$
& No & Yes & Yes\\
$R^2$, $R=\mathbb Z/8\mathbb Z$
& No & Yes & No\\
$R^2$, $\mathfrak m^2=0$,
$\dim_k\mathfrak m\geq2$
& No & No & No\\
\bottomrule
\end{tabular}
}
\end{remark}

\noindent\textbf{Concluding comparison.}
The local ISKM classification describes the modules
$R^a\oplus k^b$ over arbitrary local rings.
For SKM, socle bounds yield a rank-two
characterization over Noetherian rings, with
sufficiency provided by the CEFR/K\"othe theory.
The global ISKM characterization identifies
finite products of fields through the ISKM
property of $R^2$. The examples distinguish
multiplication, SKM, ISKM, and CEFR modules.

\section*{Conflict of interest}

The authors declare that they have no conflicts of interest regarding
publication of this paper.

 \section*{Declaration of AI Use} 
The authors used ChatGPT (OpenAI) to assist with English language editing and improving the clarity and presentation of the manuscript. The authors take full responsibility for the final content.

\end{document}